\documentclass[dvips,11pt]{article}

\usepackage[latin1]{inputenc}
\usepackage{amsmath,amsfonts}
\usepackage{amsthm}
\usepackage[numbers]{natbib}
\usepackage{enumerate}
\usepackage[ps2pdf,colorlinks]{hyperref}
\usepackage{graphicx}
\usepackage{pst-node,pstricks-add,pst-text,pst-3d,pst-tree}

\newcommand{\mailto}[1]{\href{mailto:#1}{\nolinkurl{#1}}}

\newcommand{\reff}[1]{(\ref{#1})}

\newcommand{\Z}{\mathbb{Z}}

\newcommand{\R}{\mathbb{R}}

\newcommand{\E}{\operatorname{E}}
\newcommand{\pr}{\operatorname{P}}

\newcommand{\norma}[1]{||#1||_{\rm{a}}}
\newcommand{\normv}[1]{||#1||_{\rm{v}}}
\newcommand{\ipa}[1]{\langle #1 \rangle_{\rm a}}
\newcommand{\Aa}{A_{\rm a}}
\newcommand{\Ap}{A_{\rm p}}
\newcommand{\half}{\frac{1}{2}}

\theoremstyle{plain}
\newtheorem{theorem}{Theorem}[section]

\newtheorem{lemma}[theorem]{Lemma}

\theoremstyle{definition}

\newtheorem{remark}[theorem]{Remark}

\theoremstyle{plain}
\newtheorem*{theorem*}{Theorem}
\newtheorem*{corollary*}{Corollary}
\newtheorem*{conjecture*}{Conjecture}
\newtheorem*{lemma*}{Lemma}
\newtheorem*{proposition*}{Proposition}

\theoremstyle{definition}
\newtheorem*{definition*}{Definition}
\newtheorem*{example*}{Example}
\newtheorem*{remark*}{Remark}
\newtheorem*{notabene*}{Nota Bene}
\newtheorem*{exercise*}{Exercise}
\newtheorem*{question*}{Question}

\begin{document}

\title{Stability of a spatial polling system \\ with greedy myopic service}


\author{
 Lasse Leskelä\thanks{
 Postal address: Aalto University,
 PO Box 11100, 00076 Aalto, Finland.
 URL: \url{www.iki.fi/lsl/}.
 Email: \protect\mailto{lasse.leskela@iki.fi}
 }
 \and
 Falk Unger\thanks{
 Postal address: UC Berkeley,
 665 Soda Hall, CA 94720, USA.
 Email: \protect\mailto{funger@eecs.berkeley.edu}
 }
}


\date{\today}

\maketitle

\begin{abstract}
This paper studies a spatial queueing system on a circle, polled at random locations by a myopic
server that can only observe customers in a bounded neighborhood. The server operates according to
a greedy policy, always serving the nearest customer in its neighborhood, and leaving the system
unchanged at polling instants where the neighborhood is empty. This system is modeled as a
measure-valued random process, which is shown to be positive recurrent under a natural stability
condition that does not depend on the server's scan radius. When the interpolling times are
light-tailed, the stable system is shown to be geometrically ergodic. The steady-state behavior of
the system is briefly discussed using numerical simulations and a heuristic light-traffic
approximation.
\end{abstract}





\vspace{1ex}


\vspace{1ex}

\section{Introduction}
\label{sec:Introduction}

This paper studies a spatial queueing system where customers arrive to random locations in space
that are a priori unknown to the server. The server operates sequentially in time by scanning a
bounded neighborhood of a randomly chosen location, and serving the nearest customer it observes.
After a service completion or an unsuccessful scan, a new scan is performed. Such models are
motivated by applications, such as disk storage or wireless sensor networks, where the time
required for the server to find a customer may have a major impact on the customer's sojourn time
in the system.

To simplify the analysis, we will restrict to systems where service times are negligible compared
to scanning times, by assuming all service times to be identically zero. Observe that in the
opposite case where the service times dominate, the number of customers in the system behaves like
the standard single-server queue. A further simplification is to assume that the customers
arriving during a scanning period are ignored from the ongoing scan. We expect that this
assumption is valid when the scan radius is small. The resulting system can be rephrased as
follows: A server polls the system sequentially at random locations, and either serves the nearest
observed customer immediately, or leaves the system unchanged if there are no customers within the
scan radius. Instead of a single polling server, we may alternatively think that there is an
infinite stream of servers, each of which can only perform one scan during its operational
lifetime.

The service policy is greedy in the sense that the server always aims to serve the nearest
customer. We call this model \emph{greedy polling server}, to distinguish it from spatial queueing
systems where the server travels in space towards a nearest customer (\emph{greedy traveling
server}). Although a natural choice for many applications, the greedy policy is hard to analyze
mathematically, as confirmed by earlier studies on greedy traveling
servers~\cite{altman1994,bertsimas1991,coffman1987,foss1996,kroese1994,kroese1996}. Coffman and
Gilbert~\cite{coffman1987} conjectured that the greedy traveling server on a circle is stable when
the traffic intensity is less than one, regardless of the server's traveling speed. Kroese and
Schmidt~\cite{kroese1994} proved this statement for several nongreedy policies, and Foss and
Last~\cite{foss1996,foss1998} proved it for greedy traveling servers on a finite space (see also
Meester and Quant~\cite{meester1999}). Further, Altman and Levy~\cite{altman1994} analyzed the
stability of the so-called gated-greedy policy on convex spaces, and Altman and
Foss~\cite{altman1997} derived stability conditions for nongreedy randomly traveling servers on
general spaces. Despite these affirmative results for closely related systems, the stability of
the greedy traveling server on a circle still remains an open problem.

The distribution of customers in continuum spatial queueing systems with a nongreedy server is now
relatively well understood (Kroese and Schmidt~\cite{kroese1994};
Eliazar~\cite{eliazar2003,eliazar2005}; Kavitha and Altman~\cite{kavitha2009}), whereas analytical
results related to greedy traveling servers are scarce. Coffman and Gilbert~\cite{coffman1987}
showed that the greedy traveling server on a circle closely resembles a cyclic policy in heavy
traffic (see also Litvak and Adan~\cite{litvak2001} for a similar result). Bertsimas and van
Ryzin~\cite{bertsimas1991} found two lower bounds for the mean sojourn time in the system, which
are valid for all travel policies. Kroese and Schmidt~\cite{kroese1996} derived second-order
approximations for the number of customers and workload in light traffic.

The analytical challenges related to greedy traveling servers suggest that greedy polling servers
may be difficult to analyze as well. This is why we set a modest research goal for this paper,
namely, to characterize the stability of the system. This problem may be viewed as determining the
system's throughput capacity, defined as the maximal sustainable arrival rate for which the number
of customers in the system remains stable (stochastically bounded). Foss~\cite{fossNow} recently
presented an open problem, attributing it to V.~Anantharam, conjecturing that the greedy polling
server on a circle is stable if the arrival rate is less than the polling rate, regardless of the
server's scan radius.

In this paper we prove Anantharam's conjecture \cite[Section 3.2]{fossNow} by showing that the
greedy polling server on a circle is stable if and only if the polling rate exceeds the arrival
rate. The proof is based on presenting the system as a measure-valued Markov process, and
developing a novel quadratic Lyapunov functional on the space of finite counting measures for
which the measure-valued process has negative mean drift for large customer configurations.
Besides incrementing the collection of known provable facts on spatial queueing systems with
greedy service, our analytical results may be interesting in other application areas. For
instance, the greedy polling server can be viewed as a spatial birth-and-death process. Spatial
birth-and-death processes have usually been studied in the case where all individuals have a
constant death rate, see for example Ferrari, Fernández, and Garcia~\cite{ferrari2002}; and Garcia
and Kurtz~\cite{garcia2006}, who give sufficient conditions for stability in terms of the birth
rates. The greedy polling server differs from the above birth-and-death processes in that the
death rates of individuals are governed by the Voronoi tessellation generated by the customer
locations. Borovkov and Odell~\cite{borovkov2007} have recently studied a class of
spatial--temporal point processes based on Voronoi cells, where the number of individuals is
assumed constant over time. We expect that the quadratic Lyapunov functional presented in this
paper may turn out useful in studying the ergodicity of more general spatial birth-and-death
processes.

During the final writing stage of this article, we came across an interesting recent work of
Robert~\cite{robert2010}, who considers the same problem from a different point of view. Using
entirely different techniques (a stochastically monotone construction of a stationary solution), he
proves a weaker form of stability, stating that the system has a limiting distribution for which
the number of customers is finite almost surely. Our results based on Foster--Lyapunov drift
criteria allow to prove stronger forms of stability, such as positive Harris recurrence and
geometric ergodicity, depending on the tail behavior of the interpolling time distribution.

The rest of the paper is organized as follows. In Section~\ref{sec:System} we describe the system
as a measure-valued Markov process and derive formulas for its transition operators.
Section~\ref{sec:Stability} shows the positive recurrence of the system, and
Section~\ref{sec:GeometricErgodicity} is devoted to geometric ergodicity. In
Section~\ref{sec:Steady-state} we illustrate the system dynamics with numerical simulations
complemented with a heuristic approximation of the system in light traffic.
Section~\ref{sec:Conclusions} concludes the paper.

\section{System description}
\label{sec:System}

\subsection{Notation}

The server operates on the circle $S = \{ (x_1,x_2) \in \R^2: x_1^2 + x_2^2 = \ell^2 \}$ of
circumference $\ell > 0$, where the distance $d(x,y)$ between points $x$ and $y$ is defined as the
length of the shortest arc connecting them. The state space of the system is the set $M_+(S)$ of
finite counting measures $\zeta$ on $S$, so that $\zeta(B)$ denotes the number of customers located
at $B \subset S$. The elements of $M_+(S)$ will also be called configurations, and the zero
counting measure is called the empty configuration. The total number of customers in configuration
$\zeta$ is denoted by $\normv{\zeta} = \zeta(S)$; this quantity also equals the total variation of
$\zeta$. We equip the space $M_+(S)$ with the sigma-algebra generated by the maps $\zeta \mapsto
\zeta(B)$, with $B$ being a Borel set in $S$. For real functions $f$ on $S$ we denote
\[
 \int_S f(x) \, \zeta(dx) = \sum_{x \in \zeta} f(x) \, \zeta(\{x\}),
\]
where $x \in \zeta$ is shorthand for $\zeta(\{x\}) > 0$. For more background, see for instance
Daley and Vere-Jones~\cite{daley2003}.

\subsection{Population process in continuous time}
\label{sec:PopulationCT}

Customers arrive to the circle $S$ at uniformly distributed random locations. Assuming that the
arrival locations and interarrival times are independent, the spatial--temporal arrival process
can be described as a Poisson random measure on $\R_+ \times S$ with intensity measure $\lambda dt
\, m(dx)$, where $\lambda > 0$ denotes the mean arrival rate, $dt$ the Lebesgue measure on $\R_+$,
and $m(dx)$ the uniform distribution (Haar measure) on $S$. The server polls the circle
sequentially in time by scanning a neighborhood of radius $r > 0$ of a randomly chosen location,
and either serving the nearest observed customer immediately, or leaving the system unchanged if
the scanned neighborhood is empty. Hence the probability that a customer located at $x \in \zeta$
is served during a polling event is equal to
\[
 m( B_r \cap \Gamma_\zeta(x) ),
\]
where $B_r(x)$ denotes the open $r$-ball centered at $x$, and
\[
 \Gamma_\zeta(x) = \left\{y \in S: d(x,y) < d(x',y) \ \forall x' \in \zeta \right\}
\]
denotes the Voronoi cell of a point $x$ with respect to configuration $\zeta$. (If there are many
customers located at $x$, each of them is served with equal probability.) Assuming that the
interpolling times are independent and exponential, say with parameter $\mu$, and independent of
the arrival process, the customer population in the system is described by a Markovian spatial
birth-and-death process~\cite{garcia2006} with generator
\begin{multline*}
 A f(\zeta) = \lambda \int_S ( f(\zeta+\delta_x) - f(\zeta) ) \,  m(dx) \\
 + \mu \int_S (f(\zeta - \delta_x) - f(\zeta)) \, m( B_r(x) \cap \Gamma_\zeta(x) ) \, \zeta(dx),
\end{multline*}
where $\delta_x$ denotes the Dirac measure at $x$.

\subsection{Population process at polling instants}
\label{sec:PopulationDT}

The stability regions of many ordinary queueing systems are insensitive to the shape of the
service time distribution. We will show in Section~\ref{sec:Stability} that analogously, the shape
of the interpolling distribution does not affect the stability of our model. To show this claim,
we will from now on assume that the interpolling times follow a general distribution $G$ on $\R_+$
with a finite mean. To characterize stability in terms of finite mean hitting times into small
sets (see Lemma~\ref{the:PetiteSets} in Section~\ref{sec:Irreducibility}), it is sufficient to
study the discrete-time population process $W$ obtained by sampling the system at polling
instants, so that $W_t(B)$ denotes the number of customers in $B \subset S$ just after the $t$-th
polling instant, $t \in \Z_+$ (we assume that also the initial state $W_0$ is observed just after
a polling instant). Observe that if the mean hitting time of $W$ into a set is finite, then the
same is true also for the continuous-time population process, because the mean interpolling time
is assumed finite. Therefore, we will from now on only analyze the discrete-time population
process $W$.

The population process $W$ is a discrete-time Markov process in $M_+(S)$. Given an initial state
$\zeta \in M_+(S)$, we denote the one-step transition operator of $W$ by
\[
 Af(\zeta) = \E_\zeta f(W_1) = \E \left( f(W_1) \, | \, W_0 = \zeta \right),
\]
where $f$ is a bounded or positive measurable function on $M_+(S)$. The associated probability
kernel will be denoted by
\[
 P(\zeta,B) = A 1_B(\zeta) = \pr_\zeta( W_1 \in B ), \quad B \subset M_+(S).
\]

Because the polling locations are independent of the arrivals, we can decompose the transition
operator according to
\[
 A = \Aa \circ \Ap,
\]
where the arrival operator $\Aa$ and the polling operator $\Ap$ are defined as follows. The
operator $\Aa$ acts on bounded or positive measurable functions by
\begin{equation}
 \label{eq:ArrivalOperator}
 \Aa f(\zeta) =  \sum_{n \ge 0} A_0^n f(\zeta) \, G_\lambda(n),
\end{equation}
where
\[
 A_0 f(\zeta) = \int_S f(\zeta + \delta_x) \, m(dx)
\]
corresponds to adding one new customer to a uniform random location, and
\[
 G_\lambda(n) = \int_{\R_+} e^{-\lambda s} \frac{(\lambda s)^n}{n!} \, G(ds)
\]
is the probability that $n$ customers arrive during an interpolling time. The operator $\Ap$ is
defined by
\begin{equation}
 \label{eq:PollingOperator}
 \Ap f(\zeta)
 = f(x) (1 - k_r(\zeta)) + \sum_{x \in \zeta} f(\zeta - \delta_x) \, m(B_r(x) \cap \Gamma_\zeta(x)),
\end{equation}
where
\begin{equation}
 \label{eq:ScanSuccess}
 k_r(\zeta) = m(\cup_{x \in \zeta} B_r(x))
\end{equation}
is the probability that the server finds a customer during a scan targeted into configuration
$\zeta$.

\subsection{Irreducibility and aperiodicity}
\label{sec:Irreducibility}

The following result shows that the Markov process $W$ describing the customer population in the
system may become empty with probabilities uniformly bounded away from zero. As a consequence, the
process satisfies the irreducibility and aperiodicity properties summarized below (see Meyn and
Tweedie~\cite{meyn1993} for details).

\begin{lemma}
\label{the:PetiteSets}
The Markov process $W$ is $\phi$-irreducible and strongly aperiodic, where $\phi$ is the Dirac
measure on $M_+(S)$ assigning unit mass to the zero counting measure on $S$. Moreover, the level
sets of the form $C_n = \{ \zeta \in M_+(S): \normv{\zeta} \le n\}$, $n \ge 0$, are small for $W$.
\end{lemma}
\begin{proof}
Consider an initial configuration $\zeta$ with $\normv{\zeta} = k \le n$ customers. Then the
probability that the system is empty after $n$ polling instants is greater than the probability
that no customers arrive during the first $n$ interpolling times and the server finds a customer at
each of the $k$ first polling events. Because the probability of finding a customer in a nonempty
system is at least $m(B_r)$, it follows that
\[
 P^n(\zeta,\{\zeta_0\})
 \ge \epsilon^n m(B_r)^k
 \ge \epsilon^n m(B_r)^n,
\]
where $\zeta_0$ denotes the empty configuration and $\epsilon = \int_{\R_+} e^{-\lambda s} \,
G(ds)$ is the probability that no customers arrive during an interpolling time. Denoting $\mu =
\epsilon^n m(B_r)^n \phi$, we may thus conclude that
\[
 \inf_{\zeta \in C_n} P^n(\zeta, B) \ge \mu(B)
\]
for all measurable $B \subset M_+(S)$. Thus the level set $C_n$ is small~\cite[Section
5.2]{meyn1993}. Further, the inequality $P(\zeta_0,\{\zeta_0\}) \ge \epsilon$ implies that $W$ is
strongly aperiodic~\cite[Section 5.4]{meyn1993}. \qed \end{proof}

\section{Positive recurrence}
\label{sec:Stability}

This section is devoted to deriving the main stability results (Theorems~\ref{the:Harris}
and~\ref{the:Instability}), which together imply that the system is positive recurrent if and only
if the arrival rate of customers is strictly less than the polling rate, regardless of the server's
scan radius. We start in Section~\ref{sec:PopulationSize} by discussing why it is not sufficient to
analyze the mean drift of the population size, and introduce in Section~\ref{sec:Quadratic} a
quadratic functional for which the mean drift analysis works. Section~\ref{sec:Interpolation}
discusses a key interpolation inequality that is applied in Section~\ref{sec:PositiveHarris} to
show that the mean drift with respect to the quadratic functional is negative for large
configurations. Section~\ref{sec:Instability} summarizes the behavior of the unstable system.

\subsection{Mean drift with respect to population size}
\label{sec:PopulationSize}

A common method to prove the stochastic stability of a queueing system is show that the mean drift
of the system with respect to the number of customers is strictly negative for large
configurations. To see why this approach is not sufficient for the greedy polling system in this
paper, denote the number of customers in configuration $\zeta$ by $h(\zeta) = \normv{\zeta}$, and
recall that the mean drift of the system with respect to $h$ is defined by
\[
 D h(\zeta) = A h(\zeta) - h(\zeta),
\]
where $A$ is the one-step transition operator of the system. Recall the decomposition of $A = \Aa
\circ \Ap$ in Section~\ref{sec:PopulationDT}, and observe that
\[
 \Aa h(\zeta) = h(\zeta) + \lambda s_1,
\]
where $s_1 = \int s \, G(ds)$ is the mean interpolling time, and
\[
 \Ap h(\zeta) = h(\zeta) - k_r(\zeta),
\]
where $k_r$ is the probability of a successful scan defined by~\reff{eq:ScanSuccess}. As a
consequence,
\begin{align*}
  D h(\zeta) &= \lambda s_1 - \Aa k_r(\zeta).
\end{align*}

Consider a configuration $\zeta = n \delta_x$, where $n$ customers are located in a single point $x
\in S$. Then $k_r(\zeta) = m(B_r)$, so by conditioning on whether customers arrive or not during a
polling instant, we find that
\begin{align*}
 \Aa k_r(\zeta)
 &\le k_r(\zeta) G_\lambda(0) + (1 - G_\lambda(0)) \\
 &=   1 - G_\lambda(0)(1- m(B_r)),
\end{align*}
where $G_\lambda(0) = \int e^{-\lambda s} \, G(ds)$. Hence
\[
 D h(\zeta) \ge \lambda s_1 - 1 + G_\lambda(0) (1-m(B_r)).
\]
Because the right side above does not depend on $n$, we see that $D h(\zeta)$ can be strictly
positive for arbitrarily large configurations, if $\lambda s_1 > 1 - G_\lambda(0) (1-m(B_r))$.
Hence the mean drift with respect population size cannot be used to show that the system is stable
whenever $\lambda s_1 < 1$.

\subsection{Quadratic energy functional}
\label{sec:Quadratic}

Let $M(S)$ be the space of signed counting measures on $S$, that is, measures of the form $\zeta =
\sum_{i=1}^n z_i \delta_{x_i}$, where $z_i \in \Z$ and $x_i \in S$. The space $M(S)$ is a normed
vector space with the total variation norm $\normv{\zeta} = \sum_i |z_i|$, and the subspace of
finite positive counting measures $M_+(S)$ is a convex cone in $M(S)$.

Given $0 < a \le \ell/2$, define for $\zeta,\eta \in M(S)$,
\[
 \ipa{\zeta,\eta} = \int_S \int_S (a-d(x,y))_+ \, \zeta(dx) \, \eta(dy),
\]
and denote
\[
 \norma{\zeta} = \sqrt{ \ipa{\zeta,\zeta} }.
\]

\begin{lemma}
\label{the:InnerProduct}
The map $(\zeta,\eta) \mapsto \ipa{\zeta,\eta}$ is symmetric, bilinear, and positive semidefinite.
The map $\zeta \mapsto \norma{\zeta}$ is a seminorm on $M(S)$ satisfying
\begin{equation}
 \label{eq:NormBound}
 \norma{\zeta} \le \normv{\zeta}
\end{equation}
for all $\zeta \in M(S)$.
\end{lemma}
\begin{proof}
Clearly, $\ipa{\zeta,\eta}$ is bilinear and symmetric, so we only need to show positive
semidefiniteness. We shall prove this using a probabilistic argument, representing the bilinear
functional as an expectation of a random quantity. Observe that $(a - d(x,y))_+ = m(B_{a/2}(x) \cap
B_{a/2}(y))$, where $m$ denotes the uniform probability distribution on $S$. Hence by changing the
order of integration we see that
\begin{align*}
 \ipa{\zeta,\eta}
 &= \int_S \int_S m(B_{a/2}(x) \cap B_{a/2}(y)) \, \zeta(dx) \, \eta(dy) \\
 &= \int_S \int_S \int_S 1(x \in B_{a/2}(z)) 1(y \in B_{a/2}(z)) \, m(dz) \, \zeta(dx) \, \eta(dy) \\
 &= \int_S \int_S \int_S 1(z \in B_{a/2}(x)) 1(z \in B_{a/2}(y)) \, \zeta(dx) \, \eta(dy) \, m(dz).
\end{align*}
Hence
\[
 \ipa{\zeta,\eta} = \E  \zeta( B_{a/2}(U) ) \eta( B_{a/2}(U) ),
\]
where $U$ is a random variable uniformly distributed on $S$. Especially,
\begin{equation}
 \label{eq:Seminorm}
 \ipa{\zeta,\zeta} = \E \zeta( B_{a/2}(U) )^2,
\end{equation}
which shows that $\ipa{\zeta,\zeta} \ge 0$ for all $\zeta \in M(S)$. As a consequence, $\zeta
\mapsto \norma{\zeta}$ is a seminorm (see for example Rudin~\cite[4.2]{rudin1987}). The
representation~\reff{eq:Seminorm} further shows the validity of~\reff{eq:NormBound}. \qed
\end{proof}

\begin{remark}
When $a$ is chosen so that $\ell/a$ is not an integer, one could perhaps strengthen the statement
of Lemma~\ref{the:InnerProduct} by showing that $\ipa{\zeta,\eta}$ is an inner product on $M(S)$.
However, in the sequel we will only need the fact that $\norma{\zeta}$ is a seminorm.
\end{remark}

\begin{lemma}
\label{the:ClusterBall}
For any integer $n$ and any configuration $\zeta \in M_+(S)$, there exists a closed ball $B$ in $S$
with diameter $n^{-1}$ such that
\begin{equation}
\label{eq:ClusterBall}
\zeta(B) \ge n^{-1} \normv{\zeta}.
\end{equation}
\end{lemma}
\begin{proof}
Given an integer $n$, cover the unit circle with closed balls $B_1,\dots,B_n$, each having diameter
$n^{-1}$. Then
\[
 \normv{\zeta} \le \sum_{i=1}^n \zeta(B_i) \le n \max_i \zeta(B_i),
\]
which shows that $\max_i \zeta(B_i) \ge n^{-1} \normv{\zeta}$.
\qed \end{proof}

\begin{lemma}
\label{the:NormBoundPos}
For any $a>0$ for all $\zeta \in M_+(S)$,
\begin{equation}
\label{eq:NormBoundPos}
 \frac{\sqrt{a/2}}{1 + 2/a} \normv{\zeta} \le \norma{\zeta} \le \sqrt{a} \normv{\zeta}.
\end{equation}
\end{lemma}
\begin{proof}
The upper bound in~\reff{eq:NormBoundPos} follows directly by observing that $(a-d(x,y))_+ \le a$
for all $x$ and $y$. To prove the corresponding lower bound, let $n$ be an integer such that $2/a
\le n \le 2/a + 1$, and use Lemma~\ref{the:ClusterBall} to choose a closed ball $B$ with diameter
$n^{-1}$ such that~\reff{eq:ClusterBall} holds. Because $(a-d(x,y))_+ \ge a-n^{-1} \ge a/2$ for all
$x,y \in B$, it follows that
\[
 \norma{\zeta}^2 \ge \int_B \int_B (a-d(x,y))_+ \, \zeta(dx) \, \zeta(dy)
 \ge \frac{a}{2} \zeta(B)^2.
\]
Because $n \le 2/a + 1$, the lower bound now follows using~\reff{eq:ClusterBall}.
\qed \end{proof}

\subsection{Interpolation inequality}
\label{sec:Interpolation}

This section is devoted to proving a key inequality (Lemma~\ref{the:key}), which is needed for
analyzing the mean drift of the system with respect to the seminorm $||\zeta||_a$. For a point $x$
on the circle $S$ and a positive real number $a$, we denote by $x+a$ the point on the circle
obtained by traveling distance $a$ from $x$ anticlockwise on the circle, and $x-a$ the
corresponding point obtained by traveling in the clockwise direction. Moreover, we denote by
$[u,v]$ the closed arc formed by drawing a line from $u$ to $v$ moving anticlockwise on the circle,
and by $(u,v)$ the interior of $[u,v]$.

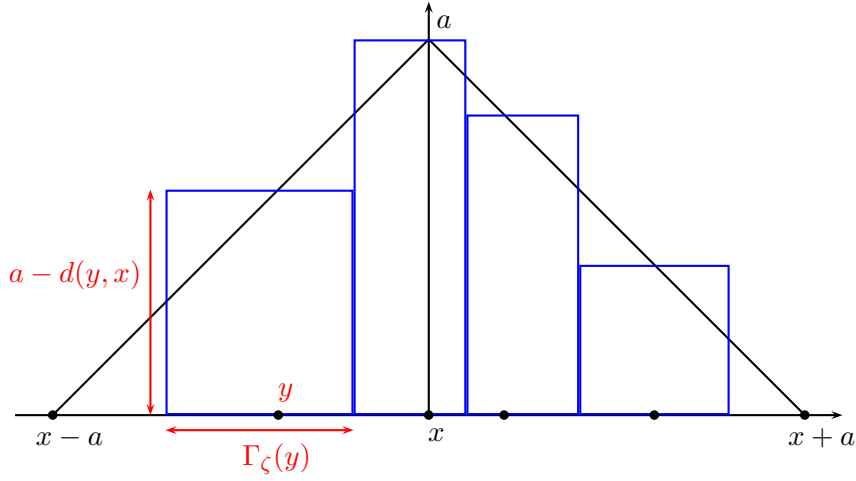
\begin{figure}[h]
\begin{center}
  \begin{pspicture}(-6,-.2)(6,6)

   \psline{->}(-5.5,0)(5.5,0)
   \psline{->}(0,0)(0,5.5)

   \psline(-5,0)(0,5)
   \psline(0,5)(5,0)

   \psframe[linestyle=solid,linecolor=blue](-3.5,0)(-1,3)
   \psframe[linestyle=solid,linecolor=blue](-1,0)(0.5,5)
   \psframe[linestyle=solid,linecolor=blue](0.5,0)(2,4)
   \psframe[linestyle=solid,linecolor=blue](2,0)(4,2)

   \psline[linecolor=red]{<->}(-3.5,-0.2)(-1,-0.2)
   \uput[-70](-2.25,-0.2){\textcolor{red}{$\Gamma_\zeta(y)$}}

   \psline[linecolor=red]{<->}(-3.7,0.0)(-3.7,3.0)
   \uput[90](-4.7,1.5){\textcolor{red}{$a-d(y,x)$}}

   \psdots[dotstyle=*, dotscale=1](-5,0)(-2,0)(0,0)(1,0)(3,0)(5,0)

   \uput[-70](-5,0){$x-a$}
   \uput[-70](0,0){$x$}
   \uput[-70](5,0){$x+a$}
   \uput[70](-2,0){\textcolor{red}{$y$}}
   \uput[50](0,5){$a$}

  \end{pspicture}
\end{center}
\caption{Nearest-neighbor interpolation of $z \mapsto (a-d(x,z))_+$. \label{fig:Triangle}}
\end{figure}

\begin{lemma}
\label{the:key2}
Let  $0 < a \le \ell/2$. Then for all $x \in S$ and for all $\zeta \in M_+(S)$ having atoms at
$x-a$, $x$, and $x+a$,
\begin{equation}
 \label{eq:KeyEq}
 \sum_{y \in \zeta} (a-d(x,y))_+ m(\Gamma_\zeta(y)) = a^2.
\end{equation}
\end{lemma}
\begin{remark}
Equation~\reff{eq:KeyEq} may interpreted in terms nearest-neighbor interpolation as follows. Given
a bounded measurable function $f$ on $S$, define its nearest-neighbor interpolant (with respect to
configuration $\zeta$) by
\[
 I_\zeta f(z) = \sum_{y \in \zeta} f(y) 1(z \in \Gamma_\zeta(y)).
\]
Then
\[
 \int_S I_\zeta f(z) \, m(dz) = \sum_{y \in \zeta} f(y) \, m(\Gamma_\zeta(y)).
\]
Given $\zeta$ and $x$ as in Lemma~\ref{the:key2}, define $f(z) = (a-d(x,z))_+$. Then the left side
in~\reff{eq:KeyEq} equals
$
 \int_S I_\zeta f(z) \, m(dz),
$ and a short calculation shows that $\int_S f(z) \, m(dz) = a^2$. Hence~\reff{eq:KeyEq} can be
rewritten as
\[
 \int_S I_\zeta f(z) \, m(dz) = \int_S f(z) \, m(dz),
\]
stating that the interpolation error made in replacing $f$ by $I_\zeta f$ is zero (see
Figure~\ref{fig:Triangle}).
\end{remark}

\begin{proof}[Lemma~\ref{the:key2}]
Assume first that $x$ is the only atom of $\zeta$ in $(x-a,x+a)$. Then
\[
 m(\Gamma_{\zeta}(x)) = m([x-a/2, x+a/2]) = a,
\]
which implies~\reff{eq:KeyEq}.

Assume next that~\reff{eq:KeyEq} holds for some configuration $\zeta$ having atoms at $x-a$, $x$,
and $x+a$. We shall show that the same is true also for $\zeta' = \zeta + \delta_z$, where $z\in S$
is arbitrary. We only need to consider the case where $z \in (x-a,x+a)$ such that $z \notin \zeta$,
because otherwise the sum on the left side of~\reff{eq:KeyEq} remains unaffected. Assume without
loss of generality that $z \in (x-a,x)$, the other case being symmetric, and let $y_1,y_2 \in
[x-a,x]$ be the points of $\zeta$ lying nearest to $z$ in the anticlockwise and clockwise
direction, respectively. Then
\begin{align*}
 m(\Gamma_{\zeta'}(z))   &= \frac{1}{2} (d(y_1,z) + d(z,y_2)), \\
 m(\Gamma_{\zeta'}(y_1)) &= m(\Gamma_\zeta(y_1)) - \half d(z,y_2), \\
 m(\Gamma_{\zeta'}(y_2)) &= m(\Gamma_\zeta(y_2)) - \half d(y_1,z).
\end{align*}

Denote by $g(x,\zeta)$ the left side of~\reff{eq:KeyEq}. Because the other atoms of $\zeta$ remain
unchanged, it follows that
\begin{multline*}
 g(x,\zeta') - g(x,\zeta)
  = \half(a-d(z,x)) ( d(y_1,z) + d(z,y_2)) \\
  - \half (a-d(y_1,x)) d(z,y_2) - \half (a-d(y_2,x)) d(y_1,z).
\end{multline*}
Because $d(y_1,x) = d(y_1,z) + d(z,x)$ and $d(y_2,x) = d(z,x) - d(y_2,z)$, the terms on the right
side of the above equation cancel out each other, so we may conclude that $g(x,\zeta') =
g(x,\zeta)$. Hence~\reff{eq:KeyEq} holds also for $\zeta' = \zeta + \delta_z$, and the proof is
complete by induction.
\qed \end{proof}

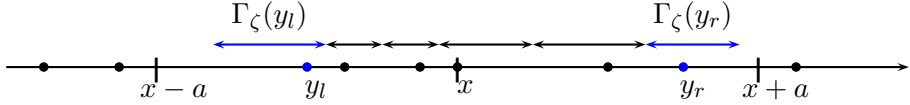
\begin{figure}[h]
\begin{center}
  \begin{pspicture}(0,-0.5)(12,1)
   \psline{->}(0,0)(12,0)
   \psdots[dotstyle=|, dotscale=2](2,0)(6,0)(10,0)
   \psdots[dotstyle=*,dotscale=1,linecolor=blue](4,0)(9,0)
   \psdots[dotstyle=*, dotscale=1](0.5,0)(1.5,0)(4.5,0)(5.5,0)(6,0)(8,0)(10.5,0)

   \uput[-70](2,0){$x-a$}
   \uput[-70](4,0){$y_l$} 
   \uput[-70](6,0){$x$}
   \uput[-70](9,0){$y_r$}
   \uput[-70](10,0){$x+a$}

   \psline[linecolor=blue]{<->}(2.75,0.3)(4.25,0.3)
   \psline{<->}(4.25,0.3)(5.00,0.3)
   \psline{<->}(5.00,0.3)(5.75,0.3)
   \psline{<->}(5.75,0.3)(7.00,0.3)
   \psline{<->}(7.00,0.3)(8.50,0.3)
   \psline[linecolor=blue]{<->}(8.5,0.3)(9.75,0.3)

   \uput[90](3.500,0.3){$\Gamma_\zeta(y_l)$}
   \uput[90](9.125,0.3){$\Gamma_\zeta(y_r)$}
  \end{pspicture}
\end{center}
\caption{Voronoi cells of $y_l$ and $y_r$. \label{fig:AddPoints}}
\end{figure}

\begin{lemma}
\label{the:key}
Assume that $0 < a \le \min(\ell/2,2r)$. Then for all $\zeta \in M_+(S)$ and for all $x \in \zeta$,
\begin{equation}
 \label{eq:keyNew}
 \sum_{y \in \zeta} (a - d(x,y))_+ m( B_r(y) \cap \Gamma_\zeta(y))
 \ge a^2.
\end{equation}
\end{lemma}
\begin{proof}
Denote the left side of~\reff{eq:keyNew} by $g(x,\zeta)$. Because $a \le \frac{\ell}{2}$, the sum
on the left of~\reff{eq:keyNew} runs over the locations $y \in \zeta$ such that $y \in [x-a,x+a]$.

Let $\zeta' = \zeta + \delta_{x-a} + \delta_{x+a}$. We shall first show that
\begin{equation}
  \label{eq:AddPoints}
  g(x,\zeta) \ge g(x,\zeta').
\end{equation}
Denote by $y_l$ and $y_r$ the atoms of $\zeta$ in $[x-a,x+a]$ located nearest to $x-a$ and $x+a$,
respectively. Then (see Figure~\ref{fig:AddPoints}) $\Gamma_{\zeta'}(y) = \Gamma_\zeta(y)$ for all
$y \in \zeta \cap (y_l,y_r)$, and moreover,
\[
 \Gamma_{\zeta'}(y_l) \subset \Gamma_\zeta(y_l) \quad \text{and} \quad \Gamma_{\zeta'}(y_r) \subset \Gamma_\zeta(y_r).
\]
Because the points $y=x-a$ and $y=x+a$ contribute nothing to the sum on the left side
of~\reff{eq:keyNew}, we may conclude that~\reff{eq:AddPoints} holds.

Further, because $a \le 2r$, it follows that $\Gamma_{\zeta'}(y) \cap B_r(y) = \Gamma_{\zeta'}(y)$
for all $y \in \zeta' \cap (x-a,x+a)$. Hence Lemma~\ref{the:key2} shows that
\[
 g(x,\zeta') = \sum_{y \in \zeta'} (a-d(x,y))_+ m(\Gamma_\zeta(y)) = a^2.
\]
In light of~\reff{eq:AddPoints}, this implies the claim.
\qed \end{proof}

\subsection{Positive Harris recurrence}
\label{sec:PositiveHarris}

The following result shows that the system is stable under the natural condition $\lambda s_1 < 1$.
Recall that $s_1$ and $s_2$ denote the first and the second moments of the interpolling time
distribution, respectively. For a converse to Theorem~\ref{the:Harris}, see
Theorem~\ref{the:Instability} in Section~\ref{sec:Instability}.

\begin{theorem}
\label{the:Harris}
Assume $\lambda s_1 < 1$ and $s_2 < \infty$. Then for all values of the scan radius $r>0$ and the
circle circumference $\ell>0$, the population process $W$ is positive Harris recurrent with
stationary distribution $\pi$ such that
\[
 \int \normv{\zeta} \, \pi(d \zeta) < \infty.
\]
Moreover, for all initial states $\zeta \in M_+(S)$,
\[
 \sup_{g: |g| \le 1 + \normv{\cdot}} \left| \E_\zeta g(W_t) - \int g \, d \pi \right| \to 0
 \qquad \text{as $t \to \infty$}.
\]
\end{theorem}
\begin{remark}
Meyn and Tweedie call such a process $f$-ergodic with $f(\zeta) = 1+\normv{\zeta}$.
\end{remark}

The proof of Theorem~\ref{the:Harris} is based on the following Foster--Lyapunov bound on the mean
drift of the system with respect to the functional $h(\zeta) = \ipa{\zeta,\zeta}$.

\begin{lemma}
\label{the:DriftQuadratic}
Assume that $s_1,s_2 < \infty$, and let $0 < a \le \min(\ell/2,2r)$. Then the mean drift of the
system with respect to $h(\zeta) = \ipa{\zeta,\zeta}$ satisfies
\[
 D h(\zeta) \le - c_1 \normv{\zeta} + c_2
\]
for all $\zeta \in M_+(S)$, where
\begin{align*}
 c_1 &= 2a^2(1 - \lambda s_1), \\
 c_2 &= a(1 + \lambda s_1) + a^2(\lambda^2 s_2 - 2\lambda s_1).
\end{align*}
\end{lemma}

\begin{proof}
Recall that $Dh(\zeta) = A h(\zeta) - h(\zeta)$, where $A = \Aa \circ \Ap$, and the transition
operators $\Aa$ and $\Ap$ are defined by~\reff{eq:ArrivalOperator} and~\reff{eq:PollingOperator}.
Observe that the transition operator corresponding to $n$ arrivals satisfies
\begin{align*}
 A_0^n h(\zeta)
 &= \E h(\zeta + \sum_{i=1}^n \delta_{X_i}) \\
 &= h(\zeta) + 2  \E \ipa{\zeta, \sum_{i=1}^n \delta_{X_i}} + \E \ipa{\sum_{i=1}^n \delta_{X_i}, \sum_{j=1}^n \delta_{X_j}} \\
 &= h(\zeta) + 2 n \E \ipa{\zeta, \delta_{X_1}} + n \E \ipa{\delta_{X_1}, \delta_{X_1}} + (n^2 - n) \E \ipa{\delta_{X_1},\delta_{X_2}},
\end{align*}
where $X_i$ are independent and uniformly distributed on $S$. Because the uniform distribution $m$
on $S$ is shift-invariant,
\[
 \int_S (a-d(x,y))_+ \, m(dy) = a^2
\]
for all $x \in S$, so that
\[
 \E \ipa{\zeta, \delta_{X_1}}
 = \int_S \sum_{y \in \zeta} (a-d(x,y))_+ \zeta(\{y\}) \, m(dx)
 = a^2 \normv{\zeta},
\]
and
\[
 \E \ipa{\delta_{X_1},\delta_{X_2}}
 = \int_S \int_S (a-d(x,y))_+ \, m(dx) \, m(dy) = a^2.
\]
Hence
\[
 A_0^n h(\zeta) = h(\zeta) + 2n a^2 \normv{\zeta} + n a + (n^2-n) a^2.
\]
Because
\[
 \int_{\R_+} \sum_{n\ge 0} (n^2 - n) \, e^{-\lambda s} \frac{(\lambda s)^n}{n!} \, G(ds)
 = \lambda^2 \int s^2 \, G(ds),
\]
we find using~\reff{eq:ArrivalOperator} that
\[
 \Aa h(\zeta) = h(\zeta) + 2 \lambda s_1 a^2 \normv{\zeta} + \lambda s_1 a + \lambda^2 s_2 a^2,
\]
where $s_k = \int s^k \, G(ds)$ denotes the $k$-th moment of the interpolling distribution.

To calculate $\Ap h$, observe that
\[
 h(\zeta - \delta_x) - h(\zeta) = - 2 \ipa{\zeta, \delta_x} + a,
\]
which shows that
\[
 \Ap h(\zeta)
 = h(\zeta) + a k_r(\zeta) - 2 \sum_{x \in \zeta} \ipa{\zeta, \delta_x} \, m(B_r(x) \cap
 \Gamma_\zeta(x)),
\]
where $k_r(\zeta)$ is the probability of a successful scan, defined by~\reff{eq:ScanSuccess}.
Lemma~\ref{the:key} implies that for all $\zeta \in M_+(S)$ and all $x \in \zeta$,
\[
 \sum_{y \in \zeta} (a-d(x,y))_+ \, m( B_r(y) \cap \Gamma_\zeta(y)) \ge a^2,
\]
so especially,
\begin{align*}
 \sum_{y \in \zeta} \ipa{\zeta, \delta_y} \, m(B_r(y) \cap \Gamma_\zeta(y))
 &\ge a^2 \normv{\zeta}.
\end{align*}
Hence
\[
 \Ap h(\zeta) \le h(\zeta) + a - 2 a^2 \normv{\zeta}.
\]
Because
\[
 \Aa h(\zeta) = h(\zeta) + 2 \lambda s_1 a^2 \normv{\zeta} + \lambda s_1 a + \lambda^2 s_2 a^2
\]
and
\[
 \Aa(\normv{\cdot})(\zeta) = \normv{\zeta} + \lambda s_1,
\]
we find that
\[
 A h(\zeta)
 \le h(\zeta) + 2 \lambda s_1 a^2 \normv{\zeta} + \lambda s_1 a + \lambda^2 s_2 a^2
 + a - 2 a^2 ( \lambda s_1 + \normv{\zeta}).
\]
\qed \end{proof}

\begin{proof}[Theorem~\ref{the:Harris}]
By Lemma~\ref{the:DriftQuadratic}, the mean drift of the system with respect to $v(\zeta) =
\norma{\zeta}^2$ for $0 < a \le \min(\ell/2,2r)$ satisfies
\[
 D v(\zeta) \le - c_1 \normv{\zeta} + c_2,
\]
for $c_1>0$ and $c_2$. Define $f(\zeta) = 1 + c \normv{\zeta}$ where $c=c_1/2$, and let $n$ be an
integer such that $n \ge (1+c_2)/(c_1/2)$. Then
\[
 D v(\zeta) \le - f(\zeta)
\]
for all $\zeta \in M_+(S)$ such that $\normv{\zeta} > n$. Moreover,
\[
 D v(\zeta) + f(\zeta) \le 1 + c_2
\]
for all $\zeta \in M_+(S)$. Lemma~\ref{the:PetiteSets} shows that the system is $\phi$-irreducible
and aperiodic, and the level set $C_n = \{\zeta \in M_+(S): \normv{\zeta} \le n\}$ is small and
thus petite. Hence the $f$-norm ergodic theorem of Meyn and Tweedie~\cite[Theorem 14.0.1]{meyn1993}
shows the claim.
\qed \end{proof}

\subsection{Instability}
\label{sec:Instability}

\begin{theorem}
\label{the:Instability}
If $\lambda s_1 \ge 1$, then the system is not positive recurrent, and if $\lambda s_1 > 1$, then
$\normv{W_t} \to \infty$ almost surely as $t \to \infty$, regardless of the initial state.
\end{theorem}
\begin{proof}
Denote the population process of the system by $W$, and let $W'$ be the population process of a
modified system with scan radius $r' = \ell/2$, so that the server finds a customer at polling
events whenever the system is nonempty. Then $\normv{W'}$ equals the number of customers in a
standard single-server queue with rate-$\lambda$ Poisson arrivals and service times distributed
according to $G$, observed just after service completions. It is not hard to see that there exists
a coupling of $W$ and $W'$ so that $\normv{W_t} \ge \normv{W'_t}$ for all $t \in \Z_+$ almost
surely, whenever $\normv{W_0} \ge \normv{W'_0}$ (see for instance Leskel\"a~\cite[Theorem
4.8]{leskela2010}). When $\lambda s_1 \ge 1$, it is well-known that the mean return time to zero of
$\normv{W'}$ is infinite~\cite{kendall1951}. The coupling then implies that the same is true for
the process $\normv{W}$, which shows that $W'$ is not positive recurrent. The second claim follows
using the same coupling, because $\normv{W'_t} \to \infty$ almost surely when $\lambda s_1 > 1$.
\qed \end{proof}

\section{Geometric ergodicity}
\label{sec:GeometricErgodicity}

The standard single-server M/G/1 queue is known to be geometrically ergodic, when the tail of the
service time distribution is light enough (Spieksma and Tweedie~\cite{spieksma1994}). An analogous
result is true for the population process $W$, as the next result shows.

\begin{theorem}
\label{the:GeometricErgodicity}
Assume that $\lambda s_1 < 1$ and the interpolling time distribution satisfies $\int e^{\theta s}
G(ds) < \infty$ for some $\theta > 0$. Then the system is geometrically ergodic in the sense that
there exist constants $\alpha >0$, $\beta>0$, and $c < \infty$ such that
\begin{equation}
  \label{eq:GeometricErgodicity}
 \sum_{t=0}^\infty e^{\alpha t} \sup_{g: |g| \le e^{\beta \norma{\cdot}}} \left| \E_\zeta g(W_t) - \int g \, d\pi \right|
 \le c e^{\beta \norma{\zeta}}
\end{equation}
for all initial states $\zeta \in M_+(S)$. Moreover, the stationary number of customers is
light-tailed in the sense that
\begin{equation}
  \label{eq:LightTails}
  \int e^{\gamma \normv{\zeta}} \, \pi(d\zeta) < \infty
\end{equation}
for some $\gamma > 0$.
\end{theorem}

Before proceeding with the proof of Theorem~\ref{the:GeometricErgodicity}, we will show that the
seminorm $\norma{\zeta}$ satisfies a similar Foster--Lyapunov bound (Lemma~\ref{the:DriftNorm}) as
the function $\ipa{\zeta,\zeta} = \norma{\zeta}^2$ in Lemma~\ref{the:DriftQuadratic}. Using the
Foster--Lyapunov bound for $\norma{\zeta}$, we will then proceed along similar lines as in the
proof of~\cite[Theorem 16.3.1]{meyn1993} (see also Borovkov and Hordijk~\cite{borovkov2004}) to
bound the mean drift of the system with respect to the function $e^{\beta \norma{\zeta}}$ for some
$\beta > 0$ (Lemma~\ref{the:DriftExp}), which is key to proving
Theorem~\ref{the:GeometricErgodicity}.

\begin{lemma}
\label{the:TaylorSqrt}
The function $x \mapsto \sqrt{1-x}$ satisfies
\[
 \sqrt{1-x} = 1 - \half x + R(x),
\]
where $|R(x)| \le 2^{-3/2} x^2 $ for all $x \in [-\half, \half]$.
\end{lemma}

\begin{proof}
Taylor's first order approximation shows that for all $x \in [-\half,\half]$, there exists $s \in
[0,1]$ such that
\[
 R(x) = \frac{1}{8} (1-sx)^{-3/2} x^2.
\]
Because $(1-sx)^{-3/2} \le 2^{3/2}$ for all $|x| \le \half$ and $s\in [0,1]$, the claim follows.
\qed \end{proof}

\begin{lemma}
\label{the:DriftNorm}
Assume that $\lambda s_1 < 1$ and $s_2 < \infty$, and let $0 < a \le \min(\ell/2,2r)$. Then there
exist $\alpha>0$, $b>0$, and an integer $n$ such that the mean drift of the system with respect to
the seminorm $v(\zeta) = \norma{\zeta}$ satisfies
\begin{equation}
 \label{eq:DriftNorm1}
 D v(\zeta) \le - \alpha
\end{equation}
for all $\zeta \in M_+(S)$ such that $\normv{\zeta} > n$, and
\begin{equation}
 \label{eq:DriftNorm2}
 D v(\zeta) \le b
\end{equation}
for all $\zeta \in M_+(S)$ such that $\normv{\zeta} \le n$.
\end{lemma}

\begin{proof}
Jensen's inequality shows that $A v \le (A v^2)^{1/2}$ holds pointwise on $M_+(S)$, so the mean
drift with respect to $v$ is bounded by
\begin{equation}
 \label{eq:DriftNormProof}
 D v
 \le (A v^2)^{1/2} - v \\
 =   ( v^2 + Dv^2 )^{1/2} - v.
\end{equation}
Because $\lambda s_1 < 1$, we see using Lemma~\ref{the:DriftQuadratic} that there exist $c>0$ such
that $Dv^2 (\zeta) \le - c \normv{\zeta}$ whenever $\normv{\zeta}$ is large enough. Because
$\norma{\zeta} \le \normv{\zeta}$ (Lemma~\ref{the:InnerProduct}), we see that
\[
 D v
 \le ( v^2 - c v)^{1/2} - v
 =   v \left( (1 - c v^{-1})^{1/2} - 1 \right)
\]
for all $\normv{\zeta}$ large enough. Lemma~\ref{the:NormBoundPos} shows that $v(\zeta) \to \infty$
as $\normv{\zeta} \to \infty$ in $M_+(S)$. Thus, for all $\zeta \in M_+(S)$ such that
$\normv{\zeta}$ is large enough, $c v^{-1} \le \half$, and using Lemma~\ref{the:TaylorSqrt} we see
that
\begin{align*}
 D v
 &\le v \left( -\half c v^{-1} + R(c v^{-1}) \right) \\
 &\le v \left( -\half c v^{-1} + 2^{-3/2} c^2 v^{-2} \right) \\
 &= -\half c + 2^{-3/2} c^2 v^{-1}.
\end{align*}
This shows the validity of~\reff{eq:DriftNorm1} for a suitable chosen $n$.

Lemma~\ref{the:DriftQuadratic} also shows that $Dv^2 \le c_2$ for all $\zeta \in M_+(S)$. Because
$v(\zeta) \le \normv{\zeta}$ (Lemma~\ref{the:InnerProduct}), inequality~\reff{eq:DriftNormProof}
shows that~\reff{eq:DriftNorm2} holds with $b=(n^2 + c_2)^{1/2}$.
\qed \end{proof}

\begin{lemma}
\label{the:DriftExp}
Assume that $\lambda s_1 < 1$ and $\int_{\R_+} e^{\theta s} \, G(ds) < \infty$ for some $\theta >
0$, and let $0 < a \le \min(\ell/2,2r)$. Then there exist $\alpha>0$, $\beta>0$, $b>0$, and an
integer $n$ such that the mean drift of the system with respect to $v_\beta(\zeta) =
\exp(\beta\norma{\zeta})$ satisfies
\begin{equation}
 \label{eq:DriftExp1}
 D v_\beta(\zeta) \le - \alpha v_\beta(\zeta)
\end{equation}
for all $\zeta \in M_+(S)$ such that $\normv{\zeta} > n$, and
\begin{equation}
 \label{eq:DriftExp2}
 D v_\beta(\zeta) + \alpha v_\beta(\zeta) \le b
\end{equation}
for all $\zeta \in M_+(S)$ such that $\normv{\zeta} \le n$.
\end{lemma}

\begin{proof}
Define $v_\beta(\zeta) = e^{\beta \norma{\zeta}}$ for some $\beta
> 0$, to be chosen later. Then the mean drift with respect to $v_\beta$ equals
\[
 D v_\beta(\zeta) = v_\beta(\zeta) \E_\zeta \left\{ e^{\beta (\norma{W_1} - \norma{W_0})} - 1 \right\}.
\]
Using a first order Taylor series approximation we see that
\[
 e^{\beta t} = 1 + \beta t + R(t),
\]
where the error term is bounded by $|R(t)| \le \frac{1}{2} \beta^2 t^2 e^{\beta |t|}$ for all $t
\in \R$. Because $\half t^2 \le s^{-2} e^{s|t|}$ for all $t$ and all $s > 0$, it follows by setting
$s=\beta^{1/3}$ that
\[
 |R(t)| \le \beta^{4/3} e^{(\beta + \beta^{1/3}) |t|}.
\]
This bound implies that
\[
 v_\beta(\zeta)^{-1} D v_\beta(\zeta)
 \le \beta D v(\zeta)
 + \beta^{4/3} \E_\zeta e^{(\beta + \beta^{1/3}) | \, \norma{W_1} - \norma{W_0)}|},
\]
where we denote $v(\zeta) = \norma{\zeta}$.

Let us next bound the exponential term. Because $\norma{\zeta}$ is a seminorm in the space of
signed counting measures on $S$ (Lemma~\ref{the:InnerProduct}), the triangle inequality shows that
\[
 | \, \norma{W_1} - \norma{W_0} | \le \norma{W_1 - W_0}.
\]
Let us write
\[
 W_1 - W_0 = \eta_a - \eta_p,
\]
where $\eta_a$ is a random counting measure describing the arrivals during an interpolling time,
and $\eta_p$ is a random counting measure describing the number of served customers during the
first polling instant ($\eta_p=0$ if the server sees no customers and $\eta_p = \delta_x$ for some
$x \in W_0 + \eta_a$ otherwise). Because $\norma{\cdot} \le \normv{\cdot}$ and $\normv{\eta_p} \le
1$, we find that
\[
 \norma{W_1 - W_0}
 \le \normv{\eta_a} + 1.
\]

Observe next that
\begin{align*}
 \E_\zeta e^{(\beta + \beta^{1/3}) \normv{\eta_a}}
 &= \int_{\R_+} \sum_{n=0}^\infty e^{(\beta + \beta^{1/3}) n} \, e^{-\lambda s} \frac{(\lambda s)^n}{n!} \, G(ds) \\
 &= \int_{\R_+} e^{\lambda(e^{(\beta + \beta^{1/3})}-1)s} \, G(ds).
\end{align*}
Choose now $\beta$ small enough such that $\lambda(e^{(\beta + \beta^{1/3})} - 1) \le \theta$ and
$e^{\beta + \beta^{1/3}} \le 2$. Then
\begin{equation}
 \label{eq:DriftExpBound}
 v_\beta(\zeta)^{-1} D v_\beta(\zeta)
 \le \beta D v(\zeta)
 + 2 \beta^{4/3} \hat G(\theta),
\end{equation}
where $\hat G(\theta) = \int e^{\theta s} \, G(ds)$. By Lemma~\ref{the:DriftNorm}, $D v(\zeta)$ is
strictly negative for $\normv{\zeta}$ large enough. Hence there exist $\alpha>0, \beta>0$ and $n$
such that~\reff{eq:DriftExp1} holds for $\normv{\zeta} > n$.

Inequality~\reff{eq:DriftExpBound} further shows that
\[
 D v_\beta(\zeta) + \alpha v_\beta(\zeta)
 \le \left( \beta D v(\zeta) + 2 \beta^{4/3} \hat G(\theta) + \alpha \right) v_\beta(\zeta).
\]
for all $\zeta \in M_+(S)$. Because $\norma{\zeta} \le \normv{\zeta}$ by
Lemma~\ref{the:InnerProduct}, we have the bound $v_\beta(\zeta) \le e^{\beta n}$ for $\normv{\zeta}
\le n$. Inequality~\reff{eq:DriftNorm2} in Lemma~\ref{the:DriftNorm} thus shows
that~\reff{eq:DriftExp2} holds for some $b$ large enough.
\qed \end{proof}

\begin{proof}[Theorem~\ref{the:GeometricErgodicity}]
By Lemma~\ref{the:PetiteSets}, the system is $\phi$-irreducible and aperiodic. Fix $0 < a \le
\min(\ell/2,2r)$, and choose $\alpha > 0$ and $\beta > 0$ as in Lemma~\ref{the:DriftExp}. By
Lemma~\ref{the:DriftExp}, the function $v_\beta(\zeta) = \exp(\beta\norma{\zeta})$ satisfies a
geometric drift condition for the level set $C_n = \{ \zeta \in M_+(S): \normv{\zeta} \le n\}$ for
some $n$ large enough. The set $C_n$ is small and thus petite by Lemma~\ref{the:PetiteSets}. Hence
by the geometric ergodic theorem of Meyn and Tweedie~\cite[Theorem 15.0.1]{meyn1993}, there exists
a finite number $c$ such that~\reff{eq:GeometricErgodicity} holds.

Observe next that~\reff{eq:GeometricErgodicity} implies
\[
 \left| g(\zeta) - \int g \, d\pi \right| \le c e^{\beta \norma{\zeta}}
\]
for all $\zeta \in M_+(S)$ and for all $g:M_+(S) \to \R$ such that $|g| \le \exp(\beta
\norma{\cdot})$. Hence, in light of Lemma~\ref{the:NormBoundPos}, inequality~\reff{eq:LightTails}
follows by applying the above inequality to the function $h(\zeta) = \exp(\gamma \normv{\zeta})$,
where $\gamma = \left( \frac{\sqrt{a/2}}{1+2/a} \right) \beta$.
\qed \end{proof}

\section{Population size in steady state}
\label{sec:Steady-state}

Having seen that the population process $W$ is positive recurrent for $\lambda s_1 < 1$ and $s_2 <
\infty$, it would be interesting to find out an explicit expression for the stationary
distribution of $W$, or at least for the stationary distribution of the number of customers in the
system. Following Kroese and Schmidt~\cite{kroese1994}, we could try the Laplace functional
approach. Denote the stationary distribution of $W$ by $\pi$, and denote the Laplace transform of
the stationary population size distribution by
\[
 L(\theta) = \int e^{-\theta \normv{\zeta}} \, \pi(d\zeta), \quad \theta > 0.
\]
Then a straightforward calculation shows that
\[
 L(\theta) = \hat G(\lambda (1-e^{-\theta})) \int e^{-\theta \normv{\zeta}} (1 + (e^\theta-1) k_r(\zeta)) \, \pi(\zeta),
\]
where $\hat G$ denotes the Laplace transform of $G$, and $k_r(\zeta)$ is defined
by~\reff{eq:ScanSuccess}. However, solving $L(\theta)$ from the above equation appears intractable
due to the nonlinearity of $k_r$.

\begin{figure}[h]
 \begin{center}
 \includegraphics[width=.49\textwidth]{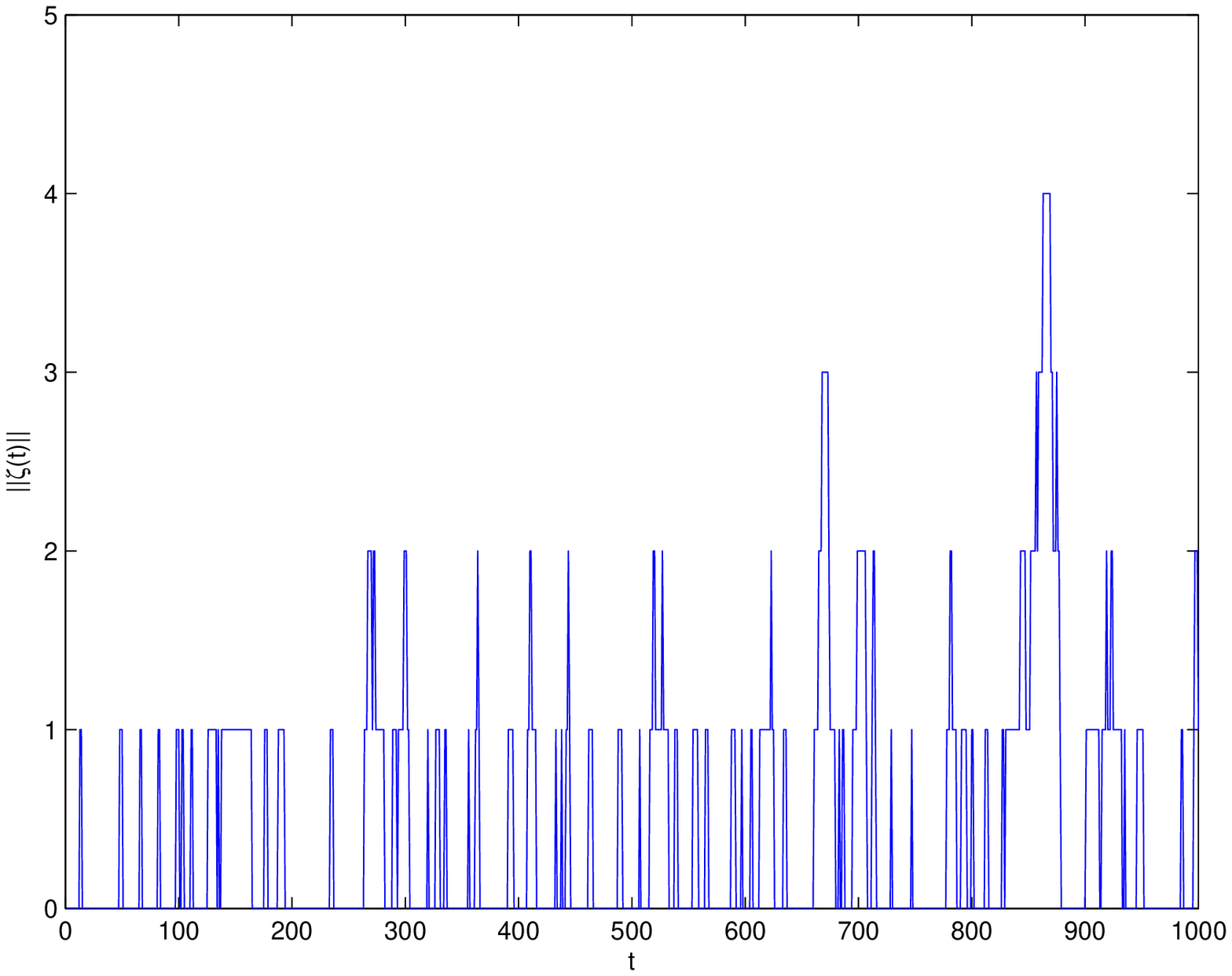}
 \includegraphics[width=.49\textwidth]{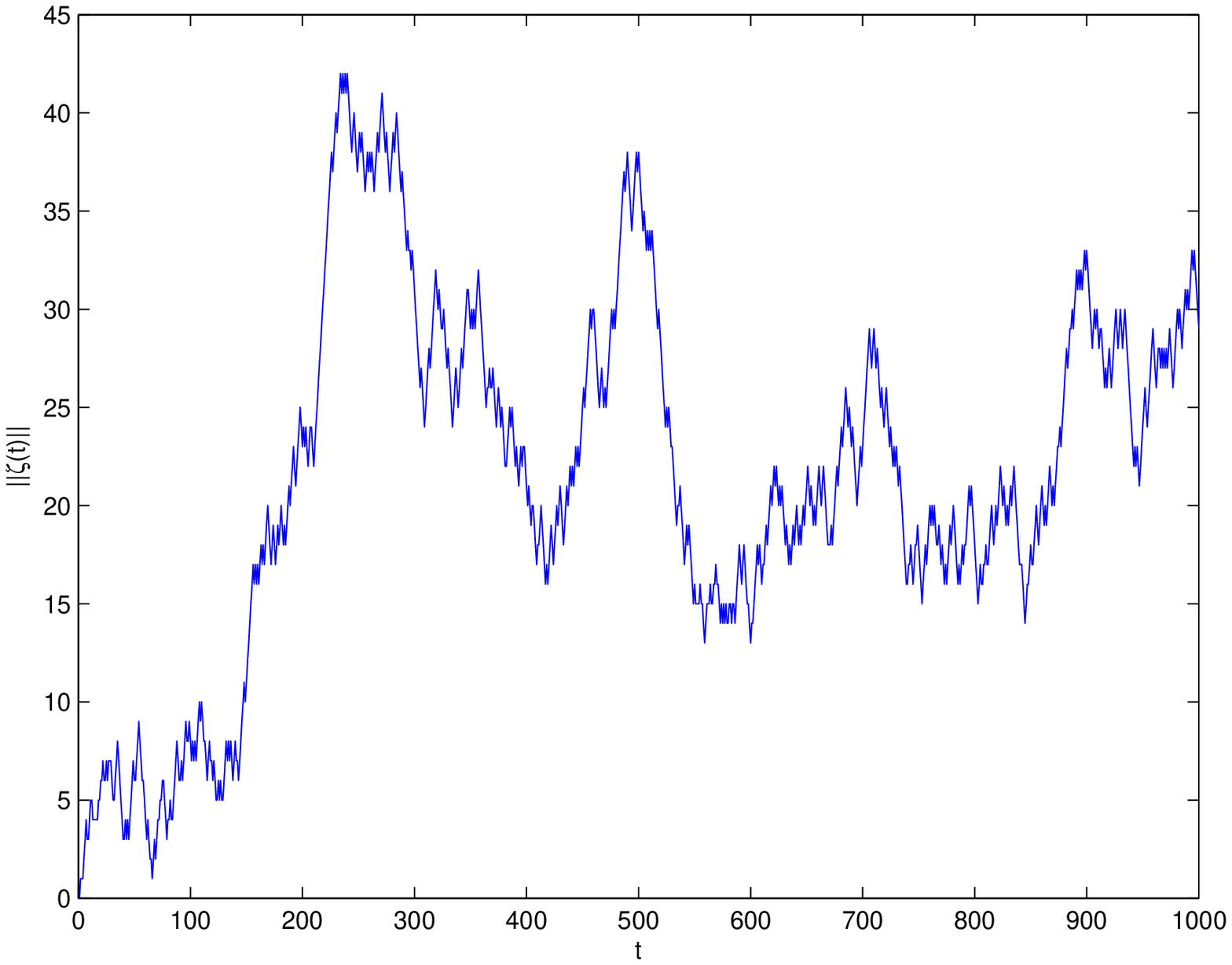}
 \caption{Population size as a function of time in light traffic ($\lambda= 0.1$, left) and heavy traffic ($\lambda= 0.9$, right).}
 \label{fig:BlueRedCircleLightTrafficPath}
 \end{center}
\end{figure}

To gain insight on the behavior of the number of customers in the system, we have numerically
simulated the system for a choice of parameter combinations with exponential interpolling times.
Figure~\ref{fig:BlueRedCircleLightTrafficPath} displays simulated paths of the population size for
$\lambda = 0.1$ (left) and $\lambda = 0.9$ (right), where $s_1 = 1$, $r = 0.1$, and $\ell = 1$.
The simulations suggest that  the system in light traffic is empty a large proportion of time,
whereas even for moderately heavy traffic ($\lambda s_1 = 0.9$), empty system appears a rare
event.

In light traffic, we can heuristically argue as follows. Assuming there is one customer present in
the system, and no new customers arrive, finding the customer requires a geometrically distributed
number of polling instants with parameter $m(B_r)$. Hence, assuming that no new customers arrive,
the time required to serve the single customer has mean $s_1/m(B_r)$. Thus we might expect that
the mean number of customers is roughly similar to the stationary probability that a renewal
on--off process with mean on-time $1/\lambda$ and mean off-time $s_1/m(B_r)$ is in on-state, so
that
\begin{equation}
  \label{eq:LightTraffic}
  \E \normv{W}
  \approx \frac{s_1/m(B_r)}{1/\lambda + s_1/m(B_r)}
  \approx \frac{\lambda s_1}{m(B_r)}.
\end{equation}

\begin{figure}[h]
 \begin{center}
 \includegraphics[width=.49\textwidth]{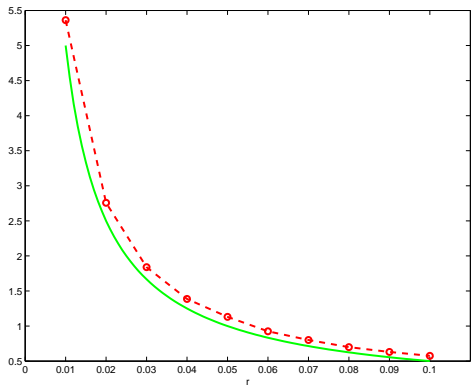}
 \caption{Simulated stationary mean population size (dashed line)
 and the light-traffic approximation (solid line) as a function of the scan radius $r$.}
 \label{fig:BlueRedCircleVaryingR}
 \end{center}
\end{figure}

To study the accuracy of the above heuristics, Figure~\ref{fig:BlueRedCircleVaryingR} plots the
approximation~\reff{eq:LightTraffic} against numerically simulated values of the stationary mean
population size for varying scan radius $r$ in a system with $\lambda=0.1$, $s_1 = 1$, and $\ell =
1$. The plot suggests that the light-traffic approximation is relatively accurate for a wide range
of scan radii.

\section{Conclusions}
\label{sec:Conclusions}
This paper studied a spatial queueing system on a circle, polled at random locations by a myopic
server that can only observe customers in a bounded neighborhood. Using a novel quadratic Lyapunov
functional of the measure-valued population process, we showed that the system is positive
recurrent under a natural stability condition, and proved the geometric ergodicity of the system
for light-tailed interpolling times. The behavior of the stationary system was discussed in terms
of numerical simulations and a heuristic light-traffic approximation. The quadratic Lyapunov
functional studied in this paper appears a promising tool for the analysis of more general spatial
birth-and-death processes.

\section*{Acknowledgments}

The results in this paper were first presented at the 15th INFORMS Applied Probability Conference,
Cornell University, 12--15 July 2009. We thank Sergey Foss for introducing this problem and for
insightful discussions during the 30th Finnish Summer School on Probability Theory in Tampere,
Finland. We also express our gratitude to two anonymous referees for their helpful remarks. This
work was initiated when L.~Leskelä was employed by the Eindhoven University of Technology and
F.~Unger by the Centrum Wiskunde \& Informatica (CWI), the Netherlands. L.~Leskelä has been
financially supported by the Academy of Finland.

\bibliographystyle{apalike}
\bibliography{lslReferences}

\end{document}